\documentclass[12pt]{article}

\usepackage{graphicx,tikz,color,url}
\usepackage{amsmath,amssymb,latexsym}
\usetikzlibrary{arrows}

\newtheorem{theorem}{Theorem}[section]

\newtheorem{lemma}[theorem]{Lemma}
\newtheorem{proposition}[theorem]{Proposition}

\newtheorem{corollary}[theorem]{Corollary}

\newcommand{\proof}{\noindent{\bf Proof.\ }}
\newcommand{\qed}{\hfill $\square$ \bigskip}

\newcommand{\smb}{S_{\rm MB}}
\newcommand{\mb}{R_{\rm MB}}

\newcommand{\diam}{{\rm diam}}

\newcommand{\cR}{{\cal R}}
\newcommand{\cS}{{\cal S}}

\newcommand{\cN}{{\cal N}}

\begin{document}
\title{Maker-Breaker resolving game played on corona products of graphs}

\author{
Tijo James $^{a,}$\thanks{Email: \texttt{tijojames@gmail.com}} 
\and 
Sandi Klav\v{z}ar $^{b,c,d,}$\thanks{Email: \texttt{sandi.klavzar@fmf.uni-lj.si}}
\and 
Dorota Kuziak $^{e,}$\thanks{Email: \texttt{dorota.kuziak@uca.es}}
\and
Savitha K S $^{f,}$\thanks{Email: \texttt{savitha@stpauls.ac.in}} 
\and 
Ambat Vijayakumar $^{g}$\thanks{Email: \texttt{vambat@gmail.com}} 
}
\maketitle

\begin{center}
$^a$ Department of Mathematics, Pavanatma College, Murickassery, India\\
\medskip

$^b$ Faculty of Mathematics and Physics, University of Ljubljana, Slovenia \\
\medskip

$^{c}$ Institute of Mathematics, Physics and Mechanics, Ljubljana, Slovenia \\
\medskip

$^d$ Faculty of Natural Sciences and Mathematics, University of Maribor, Slovenia \\
\medskip

$^{e}$ Departamento de Estad\'istica e Investigaci\'on Operativa, Universidad de C\'adiz, Algeciras, Spain \\
\medskip

$^f$ Department of Mathematics, St.Paul's  College, Kalamassery, India\\
\medskip

$^g$ Department of Mathematics, Cochin University of Science and Technology, Kochi, India
\end{center}

\newpage
\begin{abstract}
The Maker-Breaker resolving game is a game played on a graph $G$ by Resolver and Spoiler. The players taking turns alternately in which each player selects a not yet played vertex of $G$. The goal of Resolver is to select all the vertices in a resolving set of $G$, while that of Spoiler is to prevent this from happening. The outcome $o(G)$ of the game played is one of $\mathcal{R}$, $\mathcal{S}$, and $\mathcal{N}$, where $o(G)=\mathcal{R}$ (resp.\ $o(G)=\mathcal{S}$), if Resolver (resp.\ Spoiler) has a winning strategy no matter who starts the game, and $o(G)=\mathcal{N}$, if the first player has a winning strategy. In this paper, the game is investigated on corona products $G\odot H$ of graphs $G$ and $H$. It is proved that if $o(H)\in\{\mathcal{N}, \mathcal{S}\}$, then $o(G\odot H) = \mathcal{S}$. No such result is possible under the assumption $o(H) = \mathcal{R}$. It is proved that $o(G\odot P_k) = \mathcal{S}$ if $k=5$, otherwise $o(G\odot P_k) = \mathcal{R}$, and that $o(G\odot C_k) = \mathcal{S}$ if $k=3$, otherwise $o(G\odot C_k) = \mathcal{R}$. Several results are also given on corona products in which the second factor is of diameter at most $2$. 
\end{abstract}

\noindent
{\bf Keywords:} Maker-Breaker game; resolving set; Maker-Breaker resolving game; Maker-Breaker resolving number; corona product of graphs; \\

\noindent
{\bf AMS Subj.\ Class.\ (2020)}: 05C12, 05C57, 05C76

\section{Introduction}
\label{sec:intro}

Graphs in this paper are finite and simple. If 
$G=(V(G), E(G))$ is a connected   graph and $W\subseteq V(G)$, then $W$ is a {\em resolving set} of $G$ if for every pair of distinct vertices $x$ and $y$ of $G$, there exists $z\in W$ such that the shortest path distances between $x$ and $z$ and between $y$ and $z$ are different. The {\em metric dimension} $\dim(G)$ of $G$ is the minimum of the cardinalities over all resolving sets of $G$. A resolving set of cardinality $\dim(G)$ is a {\em metric basis} for $G$. This concept was independently introduced in the 1970's in~\cite{harary-1973, slater-1975} and afterwards investigated in several hundreds of papers, see the recent survey~\cite{tillquist-2023}. A fundamental reason for this incredible interest is that resolving sets and the metric dimension found numerous applications in a wide spectrum of research fields. 

The Maker-Breaker game was introduced by Erd\H{o}s and Selfridge~\cite{erdos-1973}, also in the 1970's. The game is played on an arbitrary hypergraph by two players named Maker and Breaker. The players alternately select an unplayed vertex of the hypergraph during the game. Maker's aim is to occupy all the vertices of some hyperedge, the goal of Breaker is to prevent him from doing it. The game has been extensively researched, see the book~\cite{hefetz-2014}, the recent papers~\cite{bujtas-2024, forcan-2023, hollom-2024, naimi-2023, rahman-2023}, and references therein.

In this paper we are interested in the Maker-Breaker game in which winning sets are resolving sets. This game was introduced in~\cite{kang-2021}. Closely related Maker-Breaker games with respect to distance-$k$ resolving sets and strong resolving sets were recently studied by Kang and Yi respectively in~\cite{kang-2024, kang-2024+}. In the  {\em Maker-Breaker resolving game} ({\em MBRG} for short) two players, named Resolver and Spoiler, alternately select unplayed vertices of a given graph $G$. The aim of Resolver is to select all the vertices of some resolving set of $G$, while Spoiler aims to select at least one vertex from every resolving  set of $G$. If Resolver starts the game we speak of an  {\em R-game}, otherwise we speak of an {\em S-game}. 

It was stated in~\cite{kang-2021} that for the {\em outcome} $o(G)$ of the MBRG  played on $G$ we have $o(G)\in \{\mathcal{R}, \mathcal{S}, \mathcal{N}\}$, with the following meaning: 
     \begin{itemize}
         \item $o(G)=\mathcal{R}$: Resolver has a winning strategy no matter who starts the game; 
         \item $o(G)=\mathcal{S}$: Spoiler has a winning strategy no matter who starts the game;
         \item  $o(G)=\mathcal{N}$: the first player has a winning strategy.
    \end{itemize}

Besides knowing who wins the game, it is of interest also how fast the winner can achieve this. The {\em Maker-Breaker resolving number}, $\mb(G)$, is the minimum number of moves of Resolver to win the R-game on $G$ when both players play optimally. The {\em Maker-Breaker spoiling number}, $\smb(G)$, is the minimum number of moves of Spoiler to win the R-game on $G$ when both players play optimally. For the S-game the corresponding invariants are denoted by $\mb'(G)$ and $\smb'(G)$. 

The \textit{corona product} $G\odot H$ of a connected graph $G$ and a graph $H$ is the graph obtained by taking one copy of $G$ and $n(G)$ copies of $H$ by joining the $i^{{\rm th}}$ vertex of $G$ to each vertex in the $i^{{\rm th}}$ copy of $H$, where $n(G) = |V(G)|$. For some recent studies of corona products we refer to~ \cite{alagammai-2020, dong-2021, klavzar-2022, nie-2023, nirajan-2023, verma-2024}. 

In this paper, we investigate the MBRG played on corona products. (We note in passing that the Maker-Breaker strong resolving game has already been studied on corona products in~\cite{kang-2024+}.) In the next section further definitions and  known results needed are stated. In Section~\ref{sec:outcome} we prove that if $o(H)\in\{\cN, \cS\}$, then $o(G\odot H) = \cS$. The situation when $o(H) = \cR$ is more complex. We prove that $o(G\odot P_k) = \cS$ if $k=5$, otherwise $o(G\odot P_k) = \cR$. In addition, $o(G\odot C_k) = \cS$ if $k=3$, otherwise $o(G\odot C_k) = \cR$. We also give two sufficient conditions which guarantee that $o(G\odot H) = \cR$. In the final section we consider corona products in which the second factor is of diameter at most $2$. Among other results we prove that if $\diam(H) = 2$, $o(H) = \cR$, and $\mb(H) = \mb'(H)$, then $\mb(G\odot H) = \mb'(G\odot H) = n(G)\mb(H)$. 

\section{Preliminaries}
\label{sec:prelim}

In this section we collect results to be used later on and along the way introduce some additional concepts needed. The following result will be used throughout the paper either explicitly or implicitly.

\begin{lemma} {\rm \cite[Lemma~2.2]{kang-2021}} 
\label{lem:no-skip} If the Maker-Breaker game is played on a hypergraph, then in an optimal strategy of Maker to win in the minimum number of moves it is never an advantage for him skip a move. Moreover, it never disadvantage Maker for Breaker to skip a move. 
\end{lemma}

Lemma~\ref{lem:no-skip} is known as {\em No-Skip Lemma}.

\begin{proposition} {\rm \cite[Proposition 2.3]{kang-2021}}
\label{prop:bound dimension} 
If $G$ is a connected graph, then the following hold.
\begin{enumerate}
    \item[(i)] If $o(G)=\cR$, then $\mb'(G)\geq\mb(G)\geq \dim(G)$.
    \item[(ii)] If $o(G)=\cS$, then $\smb(G)\geq\smb'(G)$.
\end{enumerate}
\end{proposition}

Let $A=\{\{u_1, v_1\},\ldots, \{u_k, v_k\}\}$ be a set of 2-subsets of $V(G)$ such that $|\bigcup_{i=1}^k{\{u_i,v_i\}}|=2k$. Then $A$ is a {\em pairing resolving set} if every set $\{x_1,\ldots, x_k\}$, where $x_i\in \{u_i,v_i\}$, is a resolving set of $G$. 

\begin{proposition} {\rm \cite[Proposition~3.3]{kang-2021}}
\label{prop:pairing} 
  If a connected graph $G$ admits a pairing resolving set, then $o(G)=\cR$.
\end{proposition}

Let $G$ be a connected graph and let $H$ be a graph. In the rest we adapt the following notation. First, $V(G) = [n(G)]$, where we used the convention that for a positive integer $k$, we write $[k] = \{1, \ldots, k\}$. Second, in $G\odot H$, the copy of $H$ which corresponds to the vertex $i\in V(G)$ will be denoted by $H_i$. We will further set $V(H_i) = \{v_1^{(i)}, \ldots, v_{n(H)}^{(i)}\}$ and denote the subgraph of $G\odot H$ induced by $V(H_i) \cup \{i\}$ by $\widehat{H}_i$. 

Let $\diam(G)$ denote the diameter of $G$, that is, the largest distance between all the pairs of vertices of $G$. The following proposition can be deduced from~\cite[Lemma~1(iv)]{yero-2011} and from the proof of~\cite[Theorem~3]{yero-2011}. 

\begin{proposition}
\label{prop:restrict-to-Hi}    
Let $G$ and $H$ be connected graphs each of order at least two. If $S$ is a resolving set of $G\odot H$, then $S\cap V(H_j)$ is a resolving set in $H_j$. Moreover, if $\diam(H) \le 2$ and $S_i$ is a resolving set of $H_i$ for $i\in [n(G)]$, then $\cup_{i=1}^{n(G)} S_i$ is a resolving set of $G\odot H$.
\end{proposition}

Let $G$ be a graph. Then a vertex $x\in V(G)$ is {\em universal} if its degree is equal to $n(G) -1$. The maximum degree of $G$ is denoted by $\Delta(G)$. For a vertex $v\in V(G)$, its open neighbourhood is $N_G(v)=\{u\in V(G):\ uv\in E(G)\}$.
$S\subseteq V(G)$ is a {\em locating set} of $G$ if $N_G(u)\cap S \ne N_G(v)\cap S$ for every two vertices $u,v\in V(G)\setminus S$. If in addition $N_G(u)\cap S \ne S$, for every $u\in V(G)\setminus S$, then $S$ is a {\em strictly locating set} of $G$. We need the following result which can be read off from~\cite[Theorem~3.6]{monsanto-2022}.

\begin{theorem}
\label{thm:locating}
Let $H$ be a connected graph and $S\subseteq V(H)$. Then $S$ is a resolving set of $K_1\odot H$ if and only if $S$ is a strictly locating set of $H$.
\end{theorem}

It is worth to observe the following consequence of Theorem~\ref{thm:locating}. 

\begin{corollary}
\label{cor:strictly-K_1}    
If Resolver has a strategy to select a strictly locating set of a connected graph $H$ in both R-game and S-game, then $o(K_1\odot H) = \cR$. 
\end{corollary}

From~\cite{kang-2021} we recall that 
if $X$ is a graph of order at least two, then $\mb(X)=1$ if and only if $X$ is a path. Similarly,  $\mb'(X)=1$ if and only if $X$ is a path. In addition, $\smb(X) \geq 2$ and $\smb'(X) \ge 2$. From these facts we can easily deduce that $\mb(G \odot H)=1= \mb'(G \odot H)$ if and only if $G\odot H \cong P_k$, where $k\in \{2,3,4\}$.

\section{Outcome of the game}
\label{sec:outcome}

In this section, we consider the outcome of the MBRG played on corona products. We first show that in the case when $o(H)\in\{\cN, \cS\}$, the outcome on $G\odot H$ is always $\cS$. The situation when $o(H) = \cR$ is more complex. This is demonstrated by the first main result of the section in which we determine the outcome for corona products $G\odot P_k$. In the second main result we determine $o(G\odot C_k)$. We conclude the section by two sufficient conditions which guarantee that $o(G\odot H) = \cR$.   

\begin{proposition}
\label{thm:outcome1} 
Let $G$ and $H$ be connected graphs  with at least two vertices. If $o(H)\in\{\cN, \cS\}$, then $o(G\odot H)=\cS$.
\end{proposition}

\proof
Assume that $o(H) \in \{\cN, \cS \}$. Then Spoiler has a winning strategy in $H$ when she starts the game. Consider an MBRG played on $G\odot H$ with Resolver as the first player. Let his first move be from the subgraph $\widehat{H}_i$. Let $j\in [n(G)]$ be an arbitrary index with $j\ne i$. Then Spoiler  responds by an optimal vertex with respect to the MBRG played on ${H}_j$. (This is possible as $n(G)\geq 2$.) In the continuation of the game, Spoiler selects vertices from ${H}_j$ according to her optimal strategy played on ${H}_j$. In view of No-Skip Lemma~\ref{lem:no-skip}, no matter, whether Resolver plays some vertices of ${H}_j$ or not, Spoiler can select a subset of vertices of ${H}_j$ such that she wins the game restricted to $\widehat{H}_j$. When the game is over, by the first assertion of Proposition~\ref{prop:restrict-to-Hi} we can conclude that Resolver cannot form a resolving set in $G\odot H$, which in turn implies that Spoiler wins the MBRG on $G\odot H$ as a second player. Therefore, by applying the No-Skip Lemma she also wins the  game as the first player. Hence $o(G\odot H)=S$. 
\qed    

If $o(H) = \cR$, then the outcome of the MBRG played on $G\odot H$ cannot be determined in general. This is demonstrated by the next theorem in which $o(G\odot P_k)$ is determined. Recall that $o(P_k) = \cR$ for $k\ge 2$ and that $o(C_k) = \cR$ for $k\ge 4$~\cite{kang-2021}. 

To prove the main results, we need the following lemma. 

\begin{lemma}
\label{lem:path-strictly-locating}   Let $\ell \ge 3$ and let $V(P_{2\ell}) = V(C_{2\ell}) = [2\ell]$. If $W=\{w_{1},\ldots, w_{\ell}\}$, where $w_{i}\in \{2i-1, 2i\}$ for $i\in [\ell]$, then $W$ is a strictly locating set of $P_{2\ell}$ as well as a a strictly locating set of $C_{2\ell}$.
\end{lemma}

\proof
Consider arbitrary vertices $x,y\in V(P_{2\ell}) \setminus W$. If $x=1$, then $2\in W$, and then, no matter whether $3\in W$ or $4\in W$, we have $N_{P_{2\ell}}(1) \cap W \ne N_{P_{2\ell}}(y)\cap W$. If $x=2$, then $1\in W$ and $1$ is clearly not in $N_{P_{2\ell}}(y)$. The cases $x=2\ell$ and $x=2\ell-1$ are symmetric. In the rest hence assume that $3\le x,y \le 2\ell - 2$.

Let $x=2j-1$. Then $2j\in W$. If $2j-2\in W$, then we easily see that $N_{P_{2\ell}}(x) \cap W \ne N_{P_{2\ell}}(y)\cap W$. And if $2j-3\in W$ then we again reach the same conclusion. The cases when $x=2j$ are done analogously. This proves the assertion for paths. 

To prove the result for cycles, note that the only obstruction for $W$ not to form a strictly locating set, is that there exist five consecutive vertices of it such that only the middle one belongs to $W$. But it readily follows that this cannot happen. 
\qed

\begin{theorem}
\label{thm:paths}    
If $G$ is a connected graph of order at least two and $k\ge 1$, then 
$$o(G\odot P_k) = \begin{cases}
\cS; & k=5, \\
\cR; & {\text otherwise}.
\end{cases}$$
\end{theorem}

\proof
We distinguish several cases. 

\medskip\noindent
{\bf Case 1}: $k=1$. \\
Recall that $v_1^{(i)}$, $i\in [n(G)]$, is  the vertex of $(P_1)_i$. Then $\{\{1,v_1 ^{(1)}\}, \ldots, \{n(G),v_1^{(n(G))}\} \}$ is a pairing resolving set, hence by Proposition~\ref{prop:pairing} we have $o(G\odot P_1) = \cR$. 

\medskip\noindent
{\bf Case 2}: $k=2$. \\
In this case the set $\{\{v_1^{(1)},v_2^{(1)}\}, \ldots, \{ v_1^{(n(G))},v_2^{(n(G))}\} \}$ is  pairing resolving, hence we can again apply Proposition~\ref{prop:pairing} to conclude that $o(G\odot P_2) = \cR$. 

\medskip\noindent
{\bf Case 3}: $k=3$. \\
In this case the strategy of Resolver is to select one of the degree two vertices in each of $(P_3)_i$. He can clearly achieve this goal and it is straightforward to see that the selected vertices form a resolving set. 

\medskip\noindent
{\bf Case 4}: $k=4$. \\
Let $v_1^{(i)}, v_2^{(i)}, v_3^{(i)}, v_4^{(i)}$ be the consecutive vertices of $(P_4)_i$. Then sets $\{v_1^{(i)},v_3^{(i)}\}$ and $\{v_2^{(i)},v_4^{(i)}\}$ form a pairing resolving set for $K_1\odot (P_4)_i$. Hence the result is true in this case also.

\medskip\noindent
{\bf Case 5}: $k=5$. \\
We are going to show that Spoiler has a winning strategy no matter who starts the game. Assume without loss of generality that the first move of Resolver is in $(\widehat{P}_5)_i$, where $i>1$. Then Spoiler replies by selecting the vertex $v_3^{(1)}$. If the next move of Resolver is $v_1^{(1)}$, then Spoiler replies by $v_5^{(1)}$. In the rest of the game Spoiler will be able to select at least one of $v_2^{(1)}$ and $v_4^{(1)}$. In either case she will win the game. In the second subcase let the next move of Resolver be  $v_2^{(1)}$. Then Spoiler replies by the vertex $v_4^{(1)}$, and afterwards she will be able to select at least one of $v_1^{(1)}$ and $v_5^{(1)}$ to win the game again. The remaining cases are symmetric, hence we have $o(G\odot P_5) = \cS$. 
 
\medskip\noindent
{\bf Case 6}: $k = 2\ell$, $\ell \ge 3$. \\
To prove that $o(G\odot P_k) = \cR$, Resolver follows the strategy that in each subgraph $(P_k)_i$ he selects one vertex from each of the sets $\{v_{2j-1}^{(i)}, v_{2j}^{(i)} \}$, $j\in [\ell]$. This strategy is indeed possible no matter who starts the game because he just follows Spoiler and as soon as she selects one vertex from $\{v_{2j-1}^{(i)}, v_{2j}^{(i)} \}$, $j\in [\ell]$, he selects the other one. 

Using the described strategy, Resolver selects in each subgraph $(P_k)_i$ a set $W_i$ which is by Lemma~\ref{lem:path-strictly-locating} a strictly locating set of $P_k \cong (P_k)_i$. By Theorem~\ref{thm:locating} we can conclude that Resolver wins the game. 

\medskip\noindent
{\bf Case 7}: $k = 2\ell+1$, $\ell \ge 3$. \\
In this case we partition the vertex set of $(P_k)_i$, $i\in [n(G)]$, as follows:  
$$Z_i = \left\{\{v_1^{(i)}, v_2^{(i)}\}, \{v_3^{(i)}, v_4^{(i)}\}, \ldots, \{ v_{2\ell-3}^{(i)}, v_{2\ell-2}^{(i)}\}, \{ v_{2\ell-1}^{(i)}, v_{2\ell}^{(i)}, v_{2\ell+1}^{(i)}\} \right\}\,.$$
The strategy of Resolver is the following. As soon as Spoiler selects a vertex from some part $\{ v_{2j-1}^{(i)}, v_{2j}^{(i)}\}$, Resolver selects the other vertex  from this part. Let now Spoiler selects for the first time a vertex from the part $\{ v_{2\ell-1}^{(i)}, v_{2\ell}^{(i)}, v_{2\ell+1}^{(i)}\}$. 

Assume that the first vertex selected from $\{ v_{2\ell-1}^{(i)}, v_{2\ell}^{(i)}, v_{2\ell+1}^{(i)}\}$ is $v_{2\ell-1}^{(i)}$. In the case that $v_{2\ell-2}^{(i)}$ has not yet been selected, Resolver selects it. Since before the end of the game Resolver will be able to select also one of $v_{2\ell}^{(i)}$ and $v_{2\ell+1}^{(i)}$, he will select in this way a strictly locating set of $(P_k)_i$. Consider next the subcase when $v_{2\ell-2}^{(i)}$ has already been selected. If $v_{2\ell-2}^{(i)}$ has been selected by Resolver (and hence $v_{2\ell-3}^{(i)}$ by Spoiler), then Resolver selects next the vertex $v_{2\ell}^{(i)}$. And if $v_{2\ell-2}^{(i)}$ has been selected by Spoiler (and hence $v_{2\ell-3}^{(i)}$ by Resolver), then Resolver selects next the vertex $v_{2\ell+1}^{(i)}$. In each of the cases the vertices selected by Resolver form a strictly locating set of $(P_k)_i$. Using Proposition~\ref{prop:pairing} once more we obtain $o(G\odot P_k) = \cR$. 

If the first vertex selected from $\{ v_{2\ell-1}^{(i)}, v_{2\ell}^{(i)}, v_{2\ell+1}^{(i)}\}$ is $v_{2\ell}^{(i)}$ or $v_{2\ell+1}^{(i)}$, then Resolver replies by selecting $v_{2\ell-1}^{(i)}$. We can then argue as above that in this way Resolver has constructed a strictly locating set of $(P_k)_i$.
\qed

The result parallel to Theorem~\ref{thm:paths} for the case $G=K_1$ is the following. 

\begin{theorem}
\label{thm:paths-by-K1}    
If $k\ge 2$, then 
$$o(K_1\odot P_k) = \begin{cases}
\cN; & k\in \{2, 5\}, \\
\cR; & {\text otherwise}.
\end{cases}$$
\end{theorem}

\proof
The proof proceeds along with the same lines as the proof of Theorem~\ref{thm:paths}, hence we skip the details here. We only emphasize that the different outcome for $P_2$ and $P_5$ comes from the fact that in $G\odot P_k$, where $n(G)\ge 2$, there are at least two different subgraphs $\widehat{H}_i$, hence in one of them Resolver can be the first player to select a vertex, while in another it is Spoiler who selects the first vertex.  
\qed

\begin{theorem}
\label{thm:cycles}    
If $G$ is a connected graph and $k\ge 3$, then 
$$o(G\odot C_k) = \begin{cases}
\cS; & k=3, \\
\cR; & {\text otherwise}.
\end{cases}$$
\end{theorem}

\proof
If $k=3$, then Spoiler can play at least two vertices in at least one subgraph $(C_3)_i$. Hence she wins the game. 

If $k\in \{4,5\}$, then Resolver follows the strategy to select two adjacent vertices in each subgraph $(C_k)_i$. Such two vertices form a strictly locating set of $(C_k)_i$ and hence in view of Theorem~\ref{thm:locating} Resolver is the winner. 

Assume next that $k\ge 6$ is even. Then the strategy of Resolver is to select in each subgraph $(C_k)_i$ a set of vertices $W_i$ as described in Lemma~\ref{lem:path-strictly-locating}. This can be achieved by accordingly responding to the moves of Spoiler. Hence we get the assertion also for even $k\ge 6$ by Theorem~\ref{thm:locating}. 

Let now $k = 2\ell + 1\ge 7$. Using the proposed notation, let $v_1^{(i)}, v_2^{(i)}, \ldots, v_{2\ell+1}^{(i)}$ be the consecutive vertices of $(C_{2\ell + 1})_i$.  We may without loss of generality assume that Spoiler starts the game by selecting the vertex $v_{2\ell+1}^{(1)}$. Resolver replies by playing the vertex $v_{2\ell}^{(1)}$. After that, Resolver continues using the  following strategy on $(C_{2\ell +1})_1$. Consider the subsets $Z_i = \{ v_{2i-1}^{(1)}, v_{2i}^{(1)} \}$, $i\in [\ell-2]$. Set also $Z' = \{v_{2\ell-3}^{(1)}, v_{2\ell-2}^{(1)}, v_{2\ell-1}^{(1)}\}$. Whenever Spoiler selects a vertex from some set $Z_i$, Resolver replies by playing the other vertex from the set. At some point Spoiler will select a vertex from $Z'$. If this vertex is $v_{2\ell-3}^{(1)}$ or $v_{2\ell-1}^{(1)}$, then Resolver replies by playing $v_{2\ell-2}^{(1)}$. And if the first vertex from $Z'$ selected by Spoiler is $v_{2\ell-2}^{(1)}$, then Resolver replies by playing $v_{2\ell-1}^{(1)}$. 
 
We claim that using the above described strategy, Resolver constructs a strictly locating set of $(C_{2\ell +1})_1$. To show it, we need to demonstrate that no four consecutive vertices were selected by Spoiler, and that there are no five consecutive vertices such that only the middle one was selected by Resolver. The first situation cannot happen because in every set $Z_i$ Resolver selected one vertex and because in $Z'$ he has also selected one vertex. The second situation could only happen if Spoiler selects $v_{2\ell+1}^{(1)}$ and $v_{1}^{(1)}$, but then since Resolver selected one of the vertices $v_{2\ell-1}^{(1)}$ and $v_{2\ell-2}^{(1)}$, this also does not happen in this case. This proves the claim.

Using the No-Skip Lemma, Resolver can select a strictly locating set in every $(C_{2\ell +1})_i$. By Theorem~\ref{thm:locating} we can conclude that Resolver is the winner. 
\qed

To conclude the section we determine two sufficient conditions which guarantee that $o(G\odot H) = \cR$. In the first we add the assumption $\diam(H) \le 2$ to $o(H)=\cR$. For the second recall that we have seen that the assumption $o(H)=\cR$ does not necessarily imply $o(G\odot H)=\cR$. 

\begin{theorem}
\label{thm:outcome-R-suffcient} 
Let $G$ and $H$ be connected graphs with  $n(G)\ge 2$ and $n(H)\geq 2$, and let at least one of the following two conditions hold: 
\begin{enumerate}
\item[(i)] $o(H)=\cR$ and $\diam(H)\leq 2$,
\item[(ii)] $o(K_1\odot H)=\cR$.
\end{enumerate}
Then $o(G\odot H)=\cR$.
\end{theorem}

\proof
(i) Since $\diam(H)\leq 2$, we infer that  $d_{H_{i}}(v_{j}^{(i)},v_{k}^{(i)})= d_{G\odot H}(v_{j}^{(i)},v_{k}^{(i)})$ holds for any two vertices $v_{j}^{(i)},v_{k}^{(i)}$ in $H_{i}$. Consider an MBRG played on $G\odot H$ with Spoiler as the first player. Spoiler selects either a vertex in $G$ say $v_{j}^{(i)}$ or a vertex in $H_i$ of $G\odot H$. Then Resolver responds optimally by choosing a non-played vertex in $H_i$ and applies his winning strategy in $H_i$ as a second player. In the continuation of the game, the strategy of Resolver is to follow Spoiler in all of the subgraphs $H_j$. 

Using the above strategy, Resolver forms a resolving set $S_i$ in each $H_i$. Since $\diam(H)\leq 2$, Proposition~\ref{prop:restrict-to-Hi} implies that  $\bigcup_i^{n(G)} S_i$ is a resolving set of $G\odot H$. Hence Resolver wins in this game as a second player. Therefore Resolver wins the game as the first player too and $o(G\odot H)=\cR$.   

(ii) Assume that $o(K_1\odot H)=\cR$. Then Resolver has a strategy to win the MBRG played on $K_1 \odot H$ no matter who starts the game. The strategy of Resolver is to apply his optimal strategy in each of the subgraphs $\widehat{H}_i$, $i\in [n(G)]$, by following Spoiler in these subgraphs. In this way Resolver forms a resolving set $S_i$ in each $\widehat{H}_i$. Since   $\bigcup_{i=1}^{n(G)}S_i$ is a resolving set of $G\odot H$, Resolver wins the game. 
\qed 

\section{Games on $G\odot H$, where $\diam(H) \le 2$}
\label{sec:with-diam-2}

In this section we consider corona products in which the second factor is of diameter at most $2$. We distinguish the cases when the first factor is $K_1$ or has at least two vertices, and bound or determine exactly the corresponding Maker-Breaker resolving numbers. 

In Theorem~\ref{thm:outcome-R-suffcient}(i) we have seen that if $G$ is a graph of order at least $2$, $H$ is a graph with $n(H)\geq 2$, $o(H)=\cR$, and $\diam(H)\leq 2$, then $o(G\odot H)=\cR$. A parallel result for the $G = K_1$ reads as follows. 

\begin{proposition}
\label{prop:diam=2-with K_1}
Let $H$ be a graph with $\diam(H) = 2$, $\Delta(H) \le n(H) - 2$, and $o(H) = \cR$. Then $o(K_1\odot H) = \cR$,  $\mb(K_1\odot H) = \mb(H)$, and $\mb'(K_1\odot H) = \mb'(H)$. \end{proposition}

\proof
In both, R-game and S-game, played on $K_1\odot H$, Resolver plays only on the subgraph $H$ by mimicking his optimal strategy in the games played on $H$. Since $\Delta(H) \le n(H) - 2$, the only universal vertex in $K_1\odot H$ is the vertex of $K_1$. Hence, having in mind that $\diam(H) = 2$, the set of vertices selected by Resolver in $H$ is also a resolving set of $K_1\odot H$. It follows that $o(K_1\odot H) = \cR$,  $\mb(K_1\odot H) \le  \mb(H)$, and $\mb'(K_1\odot H) \le  \mb'(H)$. Moreover, the equality holds in the last two inequalities because otherwise by finishing a game on $K_1\odot H$ faster, Resolver could also be able to finish a game on $H$ faster than in $\mb(H)$ moves or in  $\mb'(H)$ moves. 
\qed

The special case of Proposition~\ref{prop:diam=2-with K_1} when $H$ is the Petersen graph 
complements~\cite[Theorem~4.7]{kang-2021}. 

The assumption $\Delta(H) \le n(H) - 2$ in Proposition~\ref{prop:diam=2-with K_1} is needed to avoid the situation in which $K_1\odot H$ has more than one universal vertex. For instance, let $Y$ be the paw graph, that is, the graph obtained from $K_3$ by attaching a pendant vertex to one of its vertices. Then one can check that $o(Y) = \cR$ but $o(K_1\odot Y) = \cN$. 

A result parallel to Proposition~\ref{prop:diam=2-with K_1} for graphs $G$ of order at least two read as follows. 

\begin{proposition}
\label{prop:diam=2-with G}
Let $G$ be a connected graph with $n(G)\ge 2$, and let $H$ be a graph with $\diam(H) = 2$ and $o(H) = \cR$. Then $$\mb(G\odot H) \le \mb'(G\odot H) \le n(G)\mb'(H)\,.$$ 
If, in addition, $\mb(H) = \mb'(H)$, then 
$$\mb(G\odot H) = \mb'(G\odot H) = n(G)\mb(H)\,.$$
\end{proposition}

\proof
Let $G$ and $H$ be graphs as assumed by the proposition. By Theorem~\ref{thm:outcome-R-suffcient} we have $o(G\odot H) = \cR$. Moreover, from the proof of the same theorem we infer that Resolver will select at most $\mb'(H)$ vertices in each subgraph $\widehat{H}_i$. Since $o(G\odot H) = \cR$, we know from Proposition~\ref{prop:bound dimension}(i) that $\mb(G\odot H) \le \mb'(G\odot H)$. It follows that $\mb(G\odot H) \le \mb'(G\odot H) \le n(G)\mb'(H)$. 

Assume now that $\mb(H) = \mb'(H)$. Then Resolver selects exactly $\mb(H) = \mb'(H)$ vertices in each $\widehat{H}_i$. It follows that $\mb(G\odot H) \le \mb'(G\odot H) \le n(G)\mb(H)$. To complete the argument, we claim that $\mb(G\odot H)\ge n(G)\mb(H)$. Indeed, if this is not the case, then by using his optimal strategy, Resolver selects strictly less than $\mb(H)$ vertices in some $\widehat{H}_i$. But this would mean that Resolver can also finish the game in strictly less than $\mb(H)$ moves in $H$, which is not possible. 
 \qed

Considering the outcome of the MBRG on corona products $K_1\odot H$ is  equivalent to considering $o(G)$, where $G$ is a graph with a universal vertex. Note also that such a graph $G$ has $\diam(G)\le 2$. This problem appears very difficult in general. 

To conclude the paper we propose a stronger assumption $o(G\odot H)=\cR$ and get the following result. 

\begin{proposition}
Let $G$ and $H$ be connected graphs of order at least two and let $\diam(H)\leq 2$. If $o(G\odot H)=\cR$, then 
$$n(G) \mb(H)\leq \mb(G \odot H) \leq \mb(H)+(n(G)-1)\mb^{'}(H)\,.$$
\end{proposition}

\proof
By Proposition~\ref{thm:outcome1} it follows that $o(H)=\cR$. Therefore all MBRG parameters associated with the R-game played on $G\odot H$  and $H$ are well defined.

The lower bound follows from the fact that the union of resolving sets of $H_i$ forms a resolving set of $G\odot H$. For the upper bound, recall that $\mb(H) \leq \mb'(H)$ since $o(H) = \cR$. Since in at least one subgraph  $\widehat{H}_i$ Resolver will be the first to select a vertex, we can conclude that he will win in at most $\mb(H)+(n(G)-1)\mb'(H)$ moves.
\qed

\section*{Acknowledgements}

Sandi Klav\v zar was supported by the Slovenian Research Agency ARIS (research core funding P1-0297 and projects N1-0218, N1-0285). Dorota Kuziak was visiting the University of Ljubljana supported by “Ministerio de Educaci\'on, Cultura y Deporte”, Spain, under the “Jos\'e Castillejo” program for young researchers (reference number: CAS22/00081). A lot of work on this paper has been done during the Workshop on Games on Graphs II, June 2024, Rogla, Slovenia, the authors thank the Institute of Mathematics, Physics and Mechanics, Ljubljana, Slovenia for supporting the workshop.  

\section*{Declaration of interests}
 
The authors declare that they have no conflict of interest. 

\section*{Data availability}
 
Our manuscript has no associated data.

\end{document}